# NUMERICAL AND ANALYTICAL MODELING OF HEAT TRANSFER IN CURRENT-CARRYING CONDUCTORS USING THE HEAT EQUATION IMPLEMENTED USING FINITE-JAX

**Arturo Rodriguez[1*], Christopher Harris[2], Avinash Potluri[1], Noah L. Estrada[1], Alan M. Hernandez[1], Vyom Kumar[3], Francisco O. Aguirre Ortega[1], Vineeth Vijaya Kumar[1]**

[1]Texas A&M University - Kingsville, Kingsville, TX 78363, USA

[2]DeepVein Inc, Fremont, CA, USA

[3]Moreau Catholic High School, Hayward, CA 94544, USA

## ABSTRACT

Current-carrying conductors inevitably experience resistive heating due to the material's finite electrical conductivity. The resulting temperature distribution within the wire has essential implications for structural integrity, efficiency, and long-term reliability of electronic and power systems. In this work, we model the steady-state heat distribution in a current-carrying wire using the classical heat conduction equation. In a two-dimensional formulation, heat transport is considered both along and across the conductor. The governing partial differential equation is discretized using finite-difference methods implemented in Finite-JAX, with appropriate initial and boundary conditions, including the Dirichlet condition relevant to practical scenarios. Time integration is performed using the Euler explicit scheme, and stability constraints are systematically examined. To assess the accuracy of the numerical approach, we compare the computed temperature fields with the exact analytical solution of the heat equation for canonical geometry. Results show that the numerical prediction converges toward the analytical solution, with error norms decreasing at the expected order of accuracy. This study demonstrates how the heat equation provides a rigorous mathematical foundation for modeling resistive heating in conductors. The errors between the numerical and analytical solutions are 1.981 K, 0.8975 K, and 0.7917 K, corresponding to the L-infinity, L1, and L2 norms. For tensor computations, TPUs deliver the highest performance, surpassing GPUs and CPUs.

**KEY WORDS:** JAX, Scientific Computing, High-Performance Computing, Finite Difference, Heat Transfer

## NOMENCLATURE

| **Variables** | | | **Subscripts and Superscripts** | |
|---|---|---|---|---|
| $k$ | thermal conductivity | W/(m·K) | $x$ | x-coordinate |
| $T$ | temperature | K | $y$ | y-coordinate |
| $\dot{q}$ | volumetric heat transfer rate | W/m$^3$ | $\Omega$ | domain |
| $Q$ | source term/forcing | K/m$^2$ | $src$ | source |
| $x$ | x-coordinate | m | i | x-index |
| $y$ | y-coordinate | m | j | y-index |
| $c$ | center of circle | | n | time-index |
| $\Omega$ | domain | | num | numerical |
| $R$ | radius | m | ana | analytical |
| $\partial\Omega$ | boundary condition | | | |
| $r$ | inner radius | m | | |
| $h$ | element size | m | | |
| $L$ | length of domain | m | | |
| $E$ | error | | | |
| $N$ | nodes | | | |
| $o$ | order | | | |

## 1. INTRODUCTION

Electrical conductors experience resistance to current flow. Metal wires experience heating, a phenomenon influenced by thermal, mechanical, and operational factors that affect the performance of electronic devices and power transmission systems. Accurate temperature prediction, including how it increases and distributes within the conductor space, is critical. The properties of these conductors degrade with excessive heating, compromising structural integrity and reducing the efficiency and reliability of microelectronic circuit systems and large-scale power grids. Understanding and modeling the heating of these conductors remains a fundamental problem in applied physics and engineering practice [1].

The heat equation provides a canonical mathematical framework for describing thermal diffusion in solids. In the context of current-carrying wires, this equation governs the spatial and temporal evolution of temperature fields, leading to the generation of internal heat and the application of external boundary conditions [2, 3, 4, 5]. Analytical solutions to the heat equation exist for idealized geometries and surface conditions, providing valuable benchmarks; however, practical conductors typically involve more complex configurations that necessitate numerical approximations [6]. Finite-difference methods are widely adopted for discretization due to their simplicity, flexibility, and suitability for high-performance computing environments [7, 8, 9, 10].

Recent advances in machine learning have enabled the integration of numerical solvers with modern automatic differentiation tools. Finite-JAX leverages JAX's computational power and differentiability to build efficient, scalable solvers, extending the concept of finite differences [11, 12, 13, 14, 15, 16, 17]. This capability not only accelerates simulation but also enables sensitivity analysis and gradient-based optimization, both of which are increasingly relevant for the design and control of advanced energy and electronic systems.

In this study, we develop and verify numerical models of heat transfer in current-carrying conductors using the finite-difference method implemented with Finite-JAX. Using two-dimensional cases, we consider capturing both axial and radial temperature variations. By comparing numerical and analytical solutions, we can systematically address the influence of discretization parameters, time-stepping schemes, and surface conditions with accuracy and stability. The findings highlight the numerical framework's reliability in reproducing analytical results and provide insights into the convergence properties of the implemented schemes. More broadly, this work demonstrates the utility of differentiable programming tools in advancing the computational modeling of classical heat transfer problems. This problem is one of the development benchmarks for our Finite-JAX software.

## 2. METHODS & EQUATIONS

The governing equations of heat conduction in a current-carrying conductor are expressed as a Poisson-type problem with a localized volumetric source. Within the interior of the conductor, for $r \leq R$, the temperature field satisfies the diffusion operator augmented by a uniform volumetric heating term:

$$-k\nabla^2 T(x,y) = \dot{q}, \qquad \text{for } r \leq R, \tag{1}$$

Where $k$ denotes the thermal conductivity and $\dot{q}$ the rate of volumetric heat generation due to resistive effects, outside the conductor, the medium is taken to be isothermal, enforcing the homogeneous Dirichlet condition:

$$T(x,y) = 0, \quad \text{for } r > R \tag{2}$$

The geometry is defined in terms of a radial coordinate system centered at ($c_x$, $c_y$), with the squared radius given by:

$$r^2 = (x - c_x)^2 + (y - c_y)^2,\ R = R_{domain} \tag{3}$$

Dividing the governing equation by the thermal conductivity $k$ and introducing the normalized source term $Q = {\dot{q}}/{k}$, the heat conduction model can be cast into the canonical Poisson form. Inside the conductor domain, where $r \leq R$, the temperature satisfies:

$$\nabla^2 T(x,y) = -Q, r \leq R, T = 0, \text{for } r > R \tag{4}$$

This representation highlights the elliptic structure of the problem, where the temperature field results from a uniform volumetric source confined to the interior region, smoothly coupled to homogeneous Dirichlet boundary conditions. The formulation is equivalent to a Poisson equation with localized forcing, a prototypical problem frequently employed in numerical analysis and scientific computing to benchmark solvers. With the forcing defined piecewise as:

$$Q(x,y) = \begin{cases} {\dot{q}}/{k} & (x,y) \in \Omega \\ 0 & (x,y) \notin \Omega \end{cases} \tag{5}$$

The analytical solution can be obtained by direct integration in radial coordinates for a conductor of radius R. The solution is quadratic in the radial variable:

$$T(r) = \frac{\dot{q}}{4k}(R^2 - r^2) = \frac{Q}{4}(R^2 - r^2) \tag{6}$$

This satisfies both the governing equations in the interior and the boundary condition at $r = R$. The parabolic structure of this expression reflects the balance between internal heat generation and radial diffusion and serves as a rigorous benchmark for validating discrete approximations. The computational domain is discretized using a uniform Cartesian grid to facilitate finite-difference approximations of the governing partial differential equations. The spatial coordinates are represented as equally spaced nodal points in both the x- and y-directions:

$$x_i = i\Delta x, i = 0, \dots , N_x - 1 \tag{7}$$

$$y_j = j\Delta y, j = 0, \dots , N_y - 1 \tag{8}$$

where $N_x$ and $N_y$ denote the number of grid nodes in each direction, and $\Delta x$ and $\Delta y$ are the respective grid spacings. For a computational domain of extent *L,* the grid resolution is given by:

$$\Delta x = \Delta y = h = \frac{L}{N_x - 1} = \frac{L}{N_y - 1} \tag{9}$$

The continuous Poisson equation (Eq. 4-5) is discretized on the Cartesian mesh using a standard five-point stencil. At each interior node (*i*, *j*), the Laplacian operator is approximated by a second-order central difference, yielding:

$$-4T_{i,j} + T_{i+1,j} + T_{i-1,j} + T_{i,j+1} + T_{i,j-1} = -h^2 Q_{i,j} \tag{10}$$

This formulation (Eq. 10) defines a sparse linear system that couples each nodal temperature value to its four nearest neighbours. The discretization is formally second-order accurate in space, and due to its simplicity, serves as a canonical model for analyzing the stability and convergence of iterative solvers. To solve this system, we employ the Jacobi iterative method, which updates the temperature field at each grid point by averaging the contributions of its neighbouring nodes and the local forcing:

$$T_{i,j}^{(n+1)} = \frac{1}{4}(T_{i+1,j}^{(n)} + T_{i-1,j}^{(n)} + T_{i,j-1}^{(n)} + T_{i,j-1}^{(n)} + h^2 Q_{i,j}) \tag{11}$$

This scheme (Eq. 11) is explicit, highly parallelizable, and provides a natural implementation pathway on GPUs or other high-performance architectures. The boundary condition is enforced directly through:

$$T_{i,j}^{(n+1)} \leftarrow 0 \text{ for } (x_i, y_i) \notin \Omega \tag{12}$$

Ensuring consistency with the prescribed homogeneous Dirichlet conditions (Eq. 12).

To represent the circular geometry of the conductor within the Cartesian mesh, we employ characteristic functions as binary masks. This approach enforces the domain constraints algebraically, allowing the governing equations to be evaluated over the entire computational grid while preserving the actual circular boundary. The characteristic function of the domain, $\chi_\Omega(x_i, y_j)$, is defined as:

$$\chi_\Omega(x_i, y_i) = \begin{cases} 1 & (x_i, y_i) \in \Omega \\ 0 & otherwise \end{cases} \tag{13}$$

So that grid points inside the physical conductor are assigned unity, while all points outside are nullified (Eq. 13). The temperature field is then masked as:

$$T_{i,j} = \chi_\Omega(x_i, y_j)\tilde{T}_{i,j} \tag{14}$$

Where $\tilde{T}_{i,j}$ denotes the unconstrained nodal temperature (Eq.14). This ensures that only nodes inside the circular domain contribute to the solution, while those outsides are automatically suppressed.

Similarly, the source term is localized through an identical masking function,

Were $\chi_{src} \equiv \chi_{\Omega}$:

$$Q_{i,j} = \chi_{src}(x_i, y_i)\frac{\dot{q}}{k} \tag{15}$$

So that the volumetric heat generation is active only within the conductor (Eq. 15).

### 2.1 Convergence Criterion

The iterative solver (Jacobi scheme) terminates when the update between successive iterations falls below a prescribed tolerance. Specifically, convergence is measured using the maximum pointwise difference between two consecutive iterates:

$$\max_{(i,j)\in\Omega}\left|T_{i,j}^{(n+1)} - T_{i,j}^{(n)}\right| \leq tolerance \tag{16}$$

$$n \geq n_{max} \tag{17}$$

Ensuring that the solution has stabilized to within a specified threshold (Eq. 16). If this condition is not met, iterations continue until the maximum allowed iteration count is reached (Eq. 17). This dual stopping criterion ensures both accuracy control and computational tractability.

### 2.2 Error Norms – Validation

To assess numerical accuracy, the error at each grid point is defined as the difference between the numerical and analytical solutions, projected onto the physical domain via the characteristic mask $\chi_{\Omega}$:

$$E_{i,j} = \chi_{\Omega}(x_i, y_j)(T_{i,j}^{num} - T_{i,j}^{ana}) \tag{18}$$

(Eq.18) ensures that only nodes inside the circular conductor contribute to the error metric, while external nodes are excluded.

### 2.3 Grid Error Norms

To quantify global accuracy, several error norms are computed:

Maximum Norm ($L_{\infty}$):

$$\|E\|_{\infty} = \max_{i,j}\left|E_{i,j}\right| \tag{19}$$

This captures the most significant pointwise deviation (Eq. 19).

Euclidean Norm ($L_2$):

$$\|E\|_2 = \left(\frac{1}{N_x N_y}\sum_{i,j} E_{i,j}^2\right)^{1/2} \tag{20}$$

Providing a measure of the root-mean-square error across the domain (Eq. 20).

Average Absolute Error ($L_1$):

$$\|E\|_1 = \frac{1}{N_x N_y} \sum_{i,j} |E_{i,j}| \tag{21}$$

This quantifies the mean deviation from the exact solution (Eq. 21).

Finally, the scaling of these norms with grid resolution confirms the second-order accuracy of the finite-difference discretization:

$$\text{Error} = \mathcal{O}(h^2) \tag{22}$$

## 3. RESULTS AND DISCUSSION

### 3.1 Numerical Simulation – Heat Equation Solution

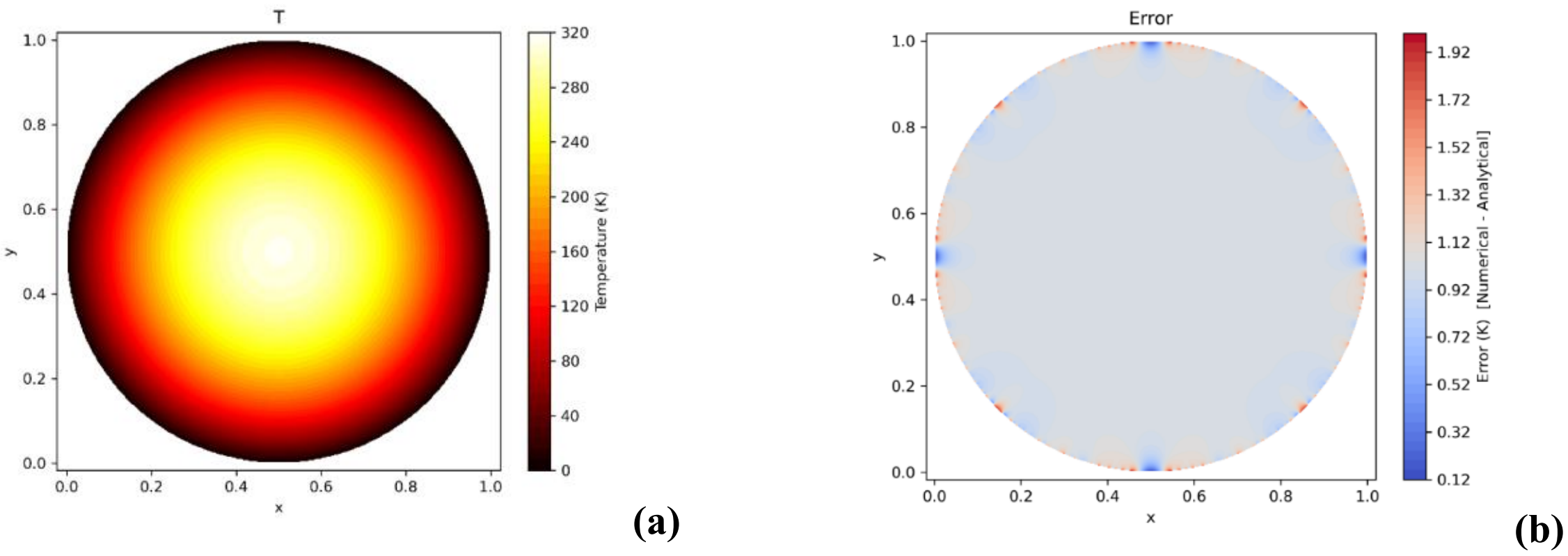

**Figs. 1** (a) Temperature Distribution and (b) Temperature Error

The figures presented illustrate the canonical solution to the two-dimensional heat equation in a circular domain. In the left figure (Figure 1), we observe the spatial temperature distribution, represented by a smooth radial profile with the highest temperature concentrated at the center and decreasing monotonically toward the surface. The numerical solution captures the expected symmetry and time-frozen conduction, reflecting the diffusivity in the separation of convection. The color map shows intense heating at the center and a smooth radial gradient, consistent with the Laplace-type analytical solution for time-frozen conduction.

On the right, the error field defines the difference between the numerical and analytical solutions. What is most notable here is the minimal magnitude of the error over most of the domain, with the discrepancies confined primarily to the right-hand bands of the circular domain. Surface effects are attributed to the discreteness and alignment of surface meshes with approximately curved geometry. The error persists at the order of machine precision in the interior, demonstrating the robustness of the discrete scheme and the spectral accuracy that we seek in high-order methods. Taken together, the results confirm the numerical solver's fidelity in reproducing the analytical driving solution, demonstrating the importance of geometric resolution on curved surfaces and minimal local errors.

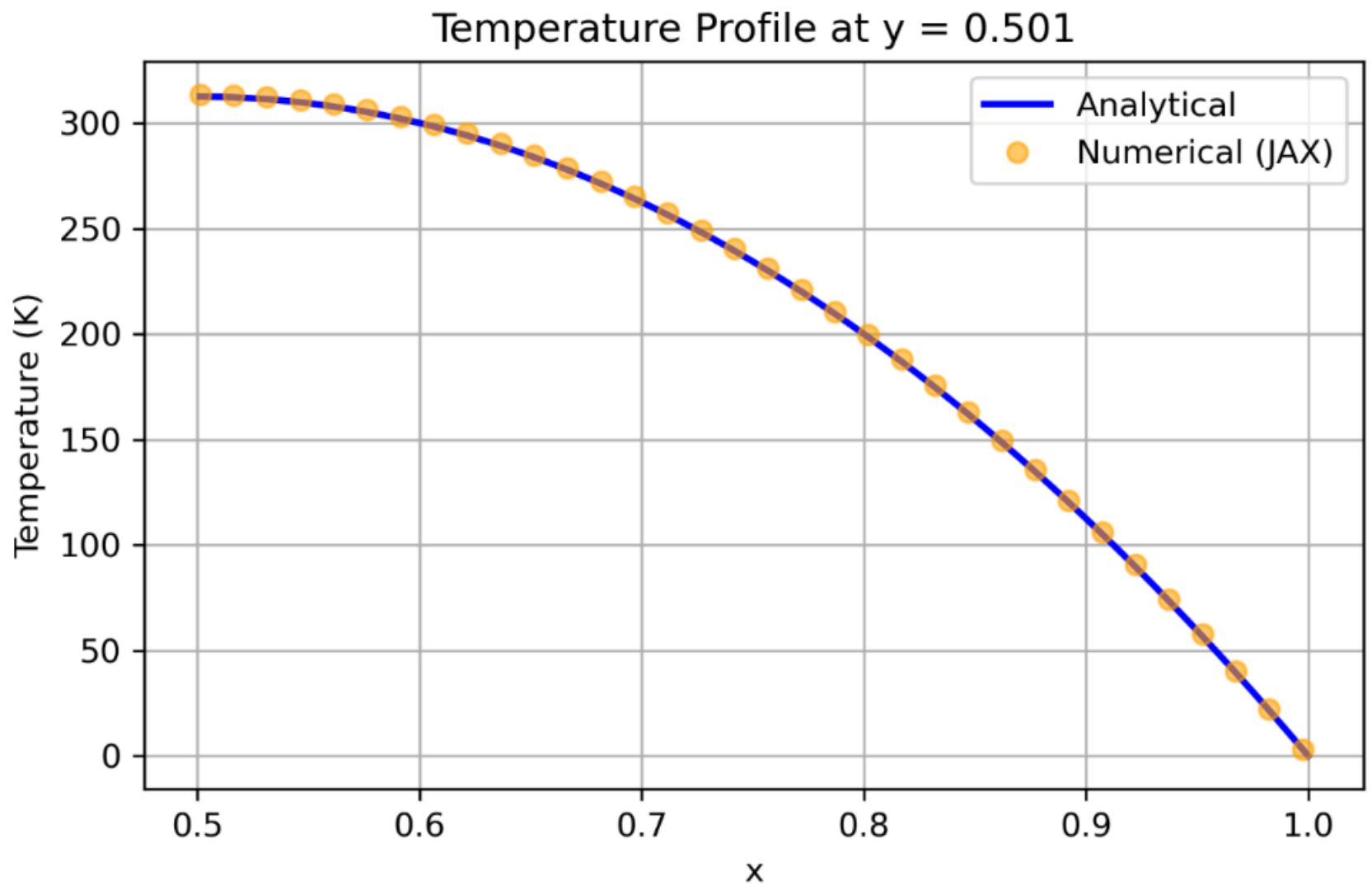

**Fig 2.** Comparison between Numerical (JAX) and Analytical Solutions

**Table 1.** Error Norms

| Maximum Norm | Euclidean Norm | Average Absolute Error |
|---|---|---|
| 1.981e+00 K | 8.975e-01 K | 7.917e-01 K |

The figure presents a one-dimensional section of the solution to the heat equation at the line y ≈ 0.5. The comparison between the analytical solution (blue line) and the numerical solution obtained using JAX (orange line) demonstrates excellent agreement across the entire domain. The profile shows a monotonic temperature drop from the interior to the surface, driven solely by diffusion.

The numerical solution is coupled with the analytical prediction solution, with minimal discrepancies. This demonstrates the accuracy and stability of the schema discretization implemented in JAX. The smoothness of the analytical curve and the superposition of the numerical points confirm the solver's consistency with the continuous model.

This validation cutoff reveals two critical aspects: (i) the solver's ability to capture the global solution structure with high fidelity, and (ii) the negligence of numerical error, which remains essential for error levels, even in regions of near-surface temperature gradients. Splicing is a vital verification step that extends the methodology to more complex geometries, time-dependent problems, and multi-physics systems.

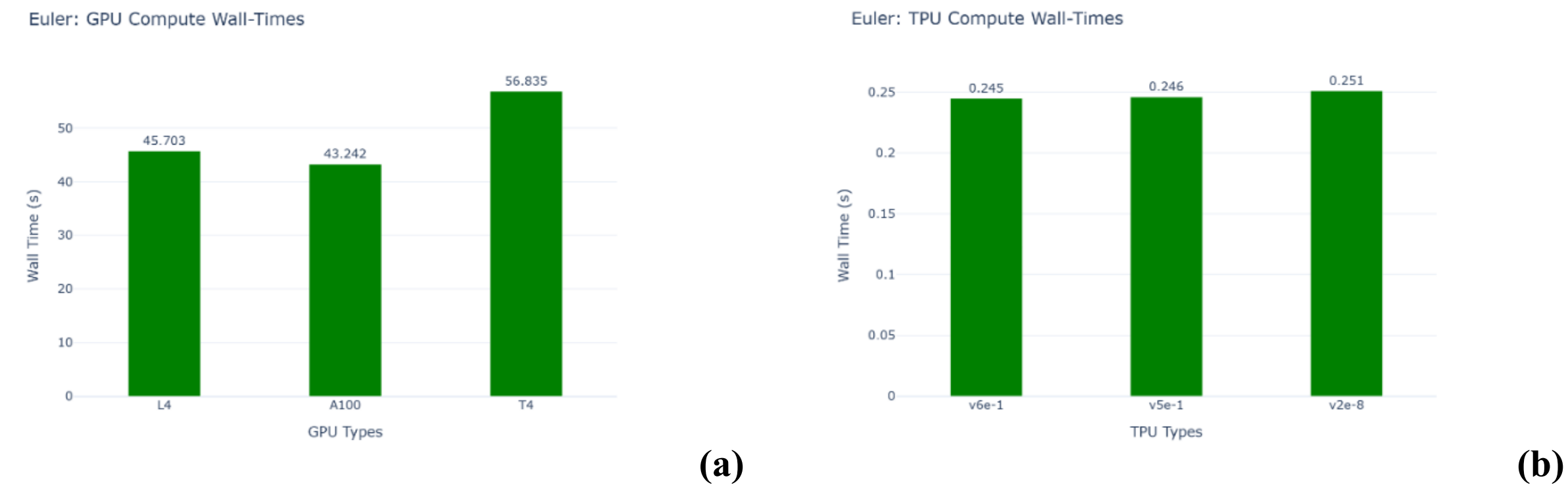

**Fig. 3** (a) Euler: GPU Wall-Times and (b) Euler: TPU Compute Wall-Times

The bar plots illustrate the clock times required to advance the Euler solver with different accelerator architectures. On the left are the GPU performances, with notable comparisons across different hardware: the NVIDIA A100 has the lowest execution time (~43 s), followed by the L4 (~46 s), and finally the T4, which has the worst performance (~57 s). This sheds light on the optimization architecture using multiple GPUs for dense linear algebra and tensor computations, which are central to partial differential equation solvers. On the right, the TPU results demonstrate the various performances, with TPU variations (v6e-1, v5e-1, and v2e-8) producing nearly identical clock times, all within 0.25 s [18, 19, 20]. This consistency across different TPUs suggests that a strong layer of the XLA compiler is effectively running the solver despite hardware differences. Taken together, the results emphasize two critical points: (i) modern GPUs perform well in partial differential equation solvers, but performance is still sensitive to architectural details, while (ii) TPUs are uniformly remarkable and scale with computational times, making them attractive for parallelizable partial differential equation workloads. These comparisons are critical for considering future directions in scientific computing and machine learning, where the trade-off among computation, portability, and energy efficiency guides hardware choice.

## 4. CONCLUSION

In this work, we have established a rigorous connection between the analytical solution of the heat equation and the numerical approximations implemented in Finite-JAX. By considering canonical conductor geometries of internal heat generation, we demonstrate discrete solutions that systematically converge to the analytical solution with scalable error norms and expected second-order rates. The use of JAX enables the integration of finite-difference solvers with differentiable programming, thereby enabling gradient-based optimization, sensitivity analysis, and coupling with data-driven models.

Beyond validating the numerical framework, our results demonstrate two implications. First, revisiting classical heat transfer problems through modern computing paradigms not only reaffirms the prediction fidelity governing partial differential equations but also reveals opportunities to accelerate GPUs and TPUs. Second, the solver's robustness against errors in geometric discretization demonstrates the importance of domain representation in preserving physical symmetry when searching for it computationally.

Looking ahead, the methodology presents the natural extension of transient formulations, combining Multiphysics and high-dimensional parameter spaces. These directions demonstrate JAX's differentiability by integrating physics-based solvers with machine learning frameworks, advancing the paradigm of scientific machine learning and computational science in thermal sciences. In this sense, the present study is not only a validation exercise but also a step toward unifying analytical theory, numerical discreteness, and modern computational methods guided by artificial intelligence in the development of next-generation heat transfer modeling.

## ACKNOWLEDGMENT

This project utilized Prof. Arturo Rodriguez's Startup funds, which Texas A&M University-Kingsville, the Department of Mechanical and Industrial Engineering, and the College of Engineering granted. This project was also supported by the DOE NNSA/MSIPP Grande CARES Consortium (GRANT13584020).

## APPENDIX

| **Simulation Parameters** | | | |
|---|---|---|---|
| ***Parameter*** | ***Value*** | ***Units*** | ***Description*** |
| k | 200.0 | W/(m-K) | Thermal Conductivity |
| L | 1.0 | m | Domain length in both x and y directions |
| Nx, Ny | 400 | | Number of grid points in the x and y directions |
| dx, dy | 0.0025 | m | Grid spacing in the x and y directions |
| cx, cy | 0.5, 0.5 | m | Coordinates of the domain center |
| R_domain | 0.5 | m | Radius of the circular computational domain |
| r0 | 0.5 | m | Radius of the internal heat source region |
| q_dot | 1x10^6 | W/m^3 | Volumetric heat generation rate |
| Q | 5000.0 | K/m^2 | Normalized heat source term used in PDE |
| dx2 | (dx)^2 | m^2 | Square of grid spacing for Jacobi update |
| tolerance | 1x10^-15 | | Convergence tolerance for iterative solver |
| max_iter | 400,000 | | Maximum number of Jacobi iterations |
| print_interval | 10,000 | | Iteration interval for debug printing |
| R | 0.5 | m | Alias for the domain radius (used in analytical solution) |
| u0 | zeros((Ny,Nx)) | K | Initial temperature field |
| u_final | computed | K | Find a numerical solution |